\documentclass[12pt, a4paper]{amsart}
\usepackage{setspace}
\usepackage[margin=2.5cm]{geometry}

\usepackage[centertags]{amsmath} 
\usepackage{amsthm}
\usepackage{amssymb}
\usepackage{mathtools} 
\usepackage{enumitem}

\usepackage{mathrsfs}
\usepackage{dutchcal}

\usepackage[english]{babel}
\usepackage[utf8x]{inputenc}

\usepackage[all]{xy}
\usepackage[bookmarks,breaklinks]{hyperref}

\usepackage[centertags]{amsmath}
\usepackage{amsthm}
\theoremstyle{plain}

\newtheorem{thm}{Theorem}[subsection]

\newtheorem{fact}[thm]{Fact} 

\newtheorem{lemma}[thm]{Lemma}
\newtheorem{prop}[thm]{Proposition}

\theoremstyle{definition}
\newtheorem{defn}[thm]{Definition}

\newcommand{\Q}{\mathbb{Q}}
\newcommand{\R}{\mathbb{R}}
\newcommand{\Z}{\mathbb{Z}}

\newcommand{\bu}{\bullet}
\renewcommand{\Im}{\operatorname{Im}}
\newcommand{\St}{\operatorname{St}}

\newcommand{\what}{\widehat}
\newcommand{\id}{\operatorname{id}}

\newcommand{\Cone}{\operatorname{Cone}}

\newcommand{\tDelta}{\tilde{\Delta}}

\newcommand{\Oo}{\operatorname{\mathcal O}}

\newcommand{\Aa}{\operatorname{\mathscr A}}
\newcommand{\Dd}{\operatorname{\mathscr D}}

\newcommand{\Cc}{\operatorname{\mathscr C}}

\newcommand{\SL}{\operatorname{SL}}

\newcommand{\ev}{\operatorname{ev}}
\newcommand{\coev}{\operatorname{coev}}

\newcommand{\gr}{\operatorname{gr}}

\newcommand{\sgn}{\operatorname{sgn}}

\newcommand{\set}[1]{\{\,#1\,\}}
\newcommand{\suchthat}[2]{\{\ #1\ \mid\ #2\ \}}

\newcommand{\Aff}{\operatorname{Aff}}

\newcommand{\Ker}{\operatorname{Ker}}

\newcommand{\wo}{\setminus}

\numberwithin{equation}{section}

\begin{document}

\bibliographystyle{alpha}

\begin{center}
  \large \textbf{Hessian metrics with distribution coefficients on a 2-sphere}\\[2ex]
\end{center}

\begin{center} 
  {\normalsize Dmitry Sustretov\footnote{ has received funding from
      the European Union's Horizon 2020 research and innovation
      programme under the Marie Sklodowska-Curie grant
      agreement No.~843100 (NALIMDIF) } \\[2ex]
  }

{\small 
\hspace{0.15\linewidth}
\begin{minipage}[t]{0.85\linewidth}
  \begin{center} {\bf Abstract}
  \end{center}  
  \footnotesize Let $\Delta$ be a 2-sphere endowed with an affine
  structure away from a finite set of points $P \subset \Delta$, and
  assume that the monodromy of the associated connection $\nabla$ on
  $\Delta \wo P$ around any point from $P$ is unipotent. I show that
  there exists a pseudo-metric tensor with distribution coefficients
  on $\Delta$ that is non-degenerate on $\Delta \wo P$ and that
  locally is of the form $\nabla d f$ for some convex function $f$. In
  particular, if $X_\infty$ is the canonical nearby fibre of a
  Type~III degeneration of K3 surfaces in Kulikov form,
  $\Delta_X \cong S^2$ is the dual intersection complex of the central
  fibre and $\Delta_X$ has simple affine structure singularities,
  existence of such ``Hessian metric'' on $\Delta_X$ implies that the
  map $H^1(\Delta_X, \Lambda^1) \to \gr^2_W H^2(X_\infty)$,
  constructed previously in \cite{sus22}, where $W$ is the monodromy
  weight filtration on $H^2(X_\infty)$ and $\Lambda^1$ is the
  push-forward of the sheaf of parallel 1-forms along the open
  embedding $\Delta \wo P \hookrightarrow \Delta$, is an isomorphism.
\end{minipage}
}
\end{center}
{
\setcounter{tocdepth}{1}
\tableofcontents

\section{Introduction}

An affine structure on a manifold $M$ is a choice of a torsion-free
flat connection $\nabla$ on its tangent bundle. This structure gives a
natural setting to define a notion of convexity. One can consider the
\emph{sheaves of $(p,q)$-forms} \cite{shima}: the sheaves of smooth
sections of the bundle $\wedge^p T^*M \otimes \wedge^q T^*M$. The
sheaves of $(p,q)$-forms form a bi-graded complex with two
differentials: $d$ and $\nabla$. A Hessian of a function $f$ is
naturally represented by a $(1,1)$-form $\nabla df$, and a smooth
function $f$ on $M$ is convex if $\nabla df$ is positive definite as a
symmetric tensor.

A notion similar to the notion of $(p,q)$-forms has received recent
interest in relation to the formalism of \emph{superforms}
\cite{lag12}, developed by Lagerberg based on an idea of Berndtsson
\cite[Section~8]{berndt}. This formalism was developed with an eye
towards application in tropical geometry and was also later adopted to
the setting of Berkovich analytic spaces, see \cite{dacl12,
  gubler16}. The superforms form a bigraded complex, similar to the
Dolbeault complex in complex geometry, with differentials $d', d''$.
The cohomology groups of the rows of this complex compute \cite{jss19}
the cohomology of certain sheaves on tropical varieties defined by
Itenberg, Katzarkov, Mikhalkin and Zharkov \cite{ikmz}. These
cohomology groups, under certain conditions, are isomorphic to even
graded pieces of the cohomology of the nearby fibres of a degeneration
of complex projective manifolds.

In \cite{sus22} I define a graded sheaf of algebras $\Lambda^\bu$ on
the dual intersection complex $\Delta_X$ of the central fibre of a
degeneration $f: X \to S$ of complex manifolds over a disc, provided
that the central fibre $f^{-1}(0)$ is a strictly normal crossings
divisor. I show that if additionally the monodromy action on a nearby
fibre $f^{-1}(t), t \neq 0$ is unipotent, and the central fibre is
projective, then for all $p,q \geq 0$ there exist morphisms
$$
H^q(\Delta_X, \Lambda^p) \to \gr^W_{2p} H^{p+q}(X, \Q).
$$
In contrast to the results of Itenberg et al. these morphisms in
general are not guaranteed to be injective or surjective. In
\cite[Section~6]{sus22} I have studied the cohomology of sheaves
$\Lambda^p$ in the case of a Type~III degeneration of K3 surfaces in
the Kulikov form, and have found a sufficient condition for
injectivity. In this case $\Delta_X$ is a 2-sphere and the sheaf
$\Lambda^1$ defines an affine structure (which was also recently
studied in \cite{eng18, aet19}) on a complement of a finite set $P$
points ($|P| \leq 24$). The morphisms above are injective subject to
existence of a certain combinatorial datum on $\Delta_X$ (a convex PL
metric on a $\Lambda^1$-torsor over $\Delta_X$). I also show that if
the monodromy of the affine structure around each point $p \in P$ is
unipotent and a primitive as an element of $\SL_2(\Z)$ (in appropriate
local coordinates), then the morphism above is surjective.

The purpose of this paper is to show that the sufficient condition for
injectivity holds using the formalism of superforms on $\Delta_X$, the
foundations for which has been laid in \cite[Sections~5.2-3]{sus22}. The
main result of this paper may also be of independent interest:

\vspace{2ex}

\textbf{Theorem~A}. Let $\Delta$ be the 2-sphere endowed with affine
structure away from finitely many points $P=\set{p_1, \ldots,
  p_n}$. Then there exists a positive symmetric $d''$-closed
$(1,1)$-supercurrent on $\Delta$ that is bounded below by a strictly
positive $(1,1)$-form.

\vspace{2ex}

For the proof of this statement see Proposition~\ref{thm-a}. By
Proposition~\ref{poincare-ddc}\ref{ddc-lemma} the current $T$ from
Theorem~A has a local potential, i.e. locally $T=d'd'' f$ for some
convex (not necessarily smooth) function $f$. If $T$ is smooth on an
open set $U \subset \Delta \wo P$, then it follows from these
consideratiosn that it defines a Hessian metric on $U$, which explains
the title of this paper.

This theorem is an analog in affine geometry of the theorem, due to
Lamari \cite{lamari} and Buchsdal \cite{buch} (see \cite[IV.3]{bhpvv}
for an exposition), that states that any smooth compact complex
surface with an even Betti number is K\"ahler. The role of the latter
condition in the complex case translates via the exact sequence of
sheaves of \emph{real} vector spaces
$$
0 \to \R \to \Oo_X \xrightarrow{\Im} \mathcal{H}_X \to 0,
$$
where $\mathcal{H}_X$ is the sheaf of pluriharmonic functions, to the
fact that the morphism
$$
H^1(X, \R) \to H^1(X, \Oo_X)
$$
in the associated long exact sequence is an isomorphism, implying
further, crucially for the proof, that $H^1(X, \mathcal{H})$ injects
into $H^2(X, \R)$. A similar role is played in the affine situation by
the sequence
$$
0 \to \R \to A^1 \to \Lambda^1 \to 0,
$$
where $A^1, \Lambda^1$ are the pushforwards from the complement of the
finite sets of affine structure singularities of the sheaves of affine
functions and parallel 1-forms, respectively. When $\Delta$ is a
2-sphere, we have $H^1(\Delta, \R)=0$ and 
$$
H^1(\Delta, A^1) \to H^1(\Delta, \Lambda^1)
$$
is injective, so the cohomology of the sheaf $A^1$ can be identified
with the cohomology clases in $H^1(\Delta, \Lambda^1)$ represented by
symmetric $(1,1)$-superforms. The key step in the proof, similarly to
the proofs of Buchsdahl and Lamari, is showing that the ``affine
Bott-Chern'' cohomology $H^1(\Delta, A^1)$ is naturally
isomorphic to its dual vector space.

\subsection{Structure of the paper}

The background on superforms on simplicial complexes is recalled in
Section~\ref{sect:background}. Section~\ref{bott-chern} contains the
main preparation steps for the proof of the Theorem~A: I describe a
modification $\tDelta$ of the 2-sphere $\Delta$ endowed with an affine
structure with singularities, for which an analogue of Serre duality
holds, prove a criterion for a existence of a K\"aler supercurrent,
and show that the sheaf $A^1$ admits a certain resolution, so that a
natural injective morphism
$H^1(\tDelta, A^1) \to H^1(\tDelta, \Lambda^1)$ exists. In
Section~\ref{main-sec} I show that $H^1(\tDelta, A^1)$ is naturally
isomorphic to its dual, and use this to prove Theorem~A.

\vspace{2ex}
\noindent\textbf{Acknowledgements}. I would like to thank Grisha
Papayanov and Julien Grivaux for useful discussions.

\section{Background}
\label{sect:background}

\label{sec:dual}


\subsection{Superforms on real vector spaces}

Let $V$ be a real vector space. We denote $V^0$ and $V^1$ two copies
of $V$ and call a \emph{superform} on $V$ a differential form on $\bar
V=V^0
\oplus V^1$ that is invariant under translations along vectors from
$V^1$.

The space of superforms admits a bidegree decompostion which is
similar to the bidegree decomposition of forms in complex geometry.
Let $\pi: V^0 \oplus V^1 \to V^0$ be the natural projeciton. We call
pullbacks $\pi^*\alpha$ of $p$-forms on $V^0$ \emph{superforms of
  bidegree $(p,0)$} on $V$. There exists a natural decomposition
$T^*\bar V = T^*V^0 \oplus T^*V^1$ and a morphism
$J: T^*\bar V \to T^*\bar V$ that swaps the two dirrect summands of
the cotangent space. This morphism induces a self-map on the
superforms. Denoting $\Aa^{p,0}$ the sheaf of $(p,0)$-superforms, we
denote
$$
\Aa^{0,q} = J \cdot \Aa^{p,0}\qquad \Aa^{p,q} := \Aa^{p,0}
\otimes_{\R} \Aa^{0,q}.
$$

The de Rham differential applied to a $(p,q)$-superform yields a
$(p+1, q)$-superform. We will denote de Rham differential
$d'$. Similarly to the differential $d^c$ of complex geometry, one
introduces the second differential
$$
d'' = J \circ d' \circ J: \Aa^{p,q} \to \Aa^{p,q+1}
$$

Pick a collection of affine functions $x_1, \ldots, x_n$ on $V$. We
will denote
$$
d'x_I = d'x_{i_1} \wedge \ldots \wedge d'x_{i_p}, \qquad d''x_J =
d''x_{j_1} \wedge \ldots \wedge d''x_{j_q}
$$
for multi-indices $I=\set{i_1, \ldots, i_p}, J=\set{j_1, \ldots, j_q},
i_1 < \ldots < i_p, j_1 < \ldots < j_q$.
One easily checks that $d',d'',J$ are given by the
following formulas
$$
\begin{array}{lll}
d'\alpha & = & \sum_{|I|=p, |J|=q} \dfrac{\partial f_{IJ}}{\partial x_i} d'x_ii
               \wedge d'x_I \wedge d''x_J\\
d''\alpha & = & (-1)^p \sum_{|I|=p, |J|=q} \dfrac{\partial f_{IJ}}{\partial x_i} d'x_ii
               \wedge d'x_I \wedge d''x_j \wedge d''x_J \\
J\alpha & = & (-1)^{pq} \sum_{|I|=p, |J|=q} f_{IJ} d'x_J \wedge d''x_I\\
\end{array}
$$

Denote $\ev: TV \otimes T^* V \to \R$ the natural pairing between
tangent vectors and covectors. The coevaluation morphism $\coev: \R
\to TV \otimes T^* V$ is the unique morphism that makes the
compositinos of the following maps
$$
\begin{array}{l}
T^* V \xrightarrow{\id \otimes \coev}  T^*V \otimes (TV \otimes T^*V)
\cong (T^*V \otimes TV) \otimes T^*V \xrightarrow{\ev \otimes
  \id} T^*V\\
TV \xrightarrow{\coev \otimes \id}   (TV \otimes T^*V) \otimes TV 
\cong TV \otimes (T^* V \otimes TV) \xrightarrow{\id \otimes
  \ev} TV\\  
\end{array}
$$
identities.

\begin{defn}[Monodromy morphism]
  Let $V$ be a real vector space. For any open set $U \subset V$
  define the morphism $N: \Aa^{p,q}(U) \to \Aa^{p-1,q+1}(U)$ to be the
  composition of the morphisms
  $$
  \begin{array}{ll}
    C^\infty(U) \otimes \wedge^p T^*V \otimes \wedge^q T^*V & \to
    C^\infty(U) \otimes \wedge^p T^*V \otimes (TV \otimes T^* V)
    \otimes \wedge^q T^*V  \\
    & \to C^\infty(U) \otimes (\wedge^p T^*V \otimes TV) \otimes (T^*
    V
    \otimes \wedge^q T^*V) \\
    & \to C^\infty(U) \otimes \wedge^{p-1} T^*V \otimes \wedge^{q+1}
    T^*V\\
  \end{array}
  $$
  where the first map is given by the coevaluation map.
  
  We also define
  $\bar N = J \circ N \circ J: \Aa^{p,q}_V \to \Aa^{p+1,q-1}_V$.
\end{defn}
The definition of the morphism $N$ is due to Yifeng Liu~\cite{liu17}.

For any ordered set $I=\set{i_1, \ldots, i_k}, i_1 < \ldots < i_k$ and
any $i = i_l \in I$ we will denote
$$
\sgn(i, I) = \dfrac{i_1 \wedge \ldots \wedge i_k}{i_l \wedge i_1
  \wedge \ldots \wedge \what{i_l} \wedge \ldots i_k} =  (-1)^{l-1}.
$$

The morphisms $N$ and $\bar N$ can be written explicitly as follows:
\begin{align*}
  N(\sum f_{IJ} d'x_I \wedge d''x_J) & = \sum_{i = 1}
  (-1)^{p-1}\sgn(i,I) d'x_{I \wo \{i\}} \wedge d''x_i \wedge d''x_J\\
  \bar N(\sum f_{IJ} d'x_I \wedge d''x_J) & = \sum_{j \in J}\sgn(j,J)
  d'x_I \wedge d'x_j \wedge d''x_{J \wo \{j\}}\\
\end{align*}

\begin{lemma}
  \label{n-diff}
  \noindent
  \begin{enumerate}
  \item\label{conj} $\bar N = J \circ N \circ J$;
  \item $[d'',N] = 0, \qquad [\bar N, d']=0$.
  \end{enumerate}
\end{lemma}

\begin{proof}
  Define
  \begin{multline*}
  \tilde{J}\left(\sum_{|I|=p,|J|=q} f_{IJ} d'x_{i_1} \wedge \ldots \wedge
  d'x_{i_p} \wedge d''x_{j_1} \wedge \ldots \wedge d''x_{j_q}\right)
  =\\
  = \sum_{|I|=p,|J|=q} f_{IJ} d'x_{i_p} \wedge \ldots \wedge
  d'x_{i_1} \wedge d''x_{j_q} \wedge \ldots \wedge d''x_{j_1}.
  \end{multline*}
  Then clearly $\bar N = \tilde{J} \circ N \circ \tilde{J}$. Now it is
  easy to see that
  $$
  J = (-1)^{pq} (-1)^{p(p-1)/2}(-1)^{q(q-1)/2} \tilde{J} =
  (-1)^{(p+q)(p+q-1)/2} J
  $$
  and so the first statement follows.
  
  The first part of the second statement is \cite[Lemma~2.2]{liu17},
  and the second one follows from the first one by \ref{conj}).
\end{proof}

\subsection{Superforms on an affine simplicial space}

Although $\Delta$ is a manifold, the singularities of the affine
structure make the sheaves of parallel forms ill-behaved, as if
$\Delta$ was singular. Notably, the sheaf of parallel 2-forms on
$\Delta$ is a constant sheaf with 1-dimensional fibres, but near any
point $x \in P$ there is only a 1-dimensional space of parallel
1-forms. It follows that parallel 2-forms cannot always be represented
as wedges of parallel 1-forms even locally. In the terminology of
\cite{sus22}, she graded sheaf of parallel forms is not regular near
the singularities of affine structure.

In order to fix this we will have to leave the category of manifolds
endowed with an affine structure with singularities. Given a
simplicial complex and a certain sheaf of functions on it that can
serve as ``affine coordinates'', one can define the sheaves of
superforms. We will call such spaces \emph{affine simplicial 
  spaces}. In particular, $\Delta$ is one, but instead of working with
$\Delta$ we will work with a polyhedral affine space $\tDelta$ which
will have the desirable properties.

\begin{defn}[Affine simplicial spaces]
  We call a pair $(\Delta, A^1)$ of a topological space and a sheaf
  of continuous functions on it an \emph{affine simplicial space} if
  $\Delta$ admits a structure of a simplicial complex and $A^1$ is a
  subsheaf of the sheaf of continuous piece-wise affine functions on
  $\Delta$ that are linear on each face.
\end{defn}

Let $\Delta$ be a surface and let $\nabla$ be a flat connection that
defines an affine structure on $\Delta$. One can define the sheaves of
superforms on $\Delta$ by putting
$$
\Aa^{p,q} = \Aa^p \otimes_{\R} \Aa^q,
$$
where $\Aa^p, \Aa^q$ are sheaves of real $p$- and $q$-forms,
respectively. The differentials $d',d''$ are given by the de Rham
differential $d$ and the connection $\nabla$, respectively.

If $(\Delta, A^1)$ is an affine simplicial space then one can define
superforms on $\Delta$ following the procedure proposed in
\cite[Section~5.3]{sus22}. First, we define forms on simplicial
complexes (or more generally, polyhedral complexes) embedded into a
vector space.

\begin{defn}[Superforms on polyhedral complexes]
  \label{sf-emb}
  Let $\Sigma$ be a polyherdal complex embedded into a real vector
  space $V$. Define the following equivalence relation on the
  superforms defines on a neighbourhood of $\Sigma$ in $V$: two forms
  $\eta$ and $\eta'$ are eqivalent if for any face
  $\sigma \subset \Sigma$, the restrictions $\eta|_{\sigma}$ and
  $\eta'|_{\sigma}$ coincide. We call such an equivalence class of
  superformsa a \emph{superform on $\Sigma$}.
\end{defn}

For any point $x \in \Delta$ there exists a neighbourhood
$U$ of $x$ and a map
$$
e: U \to (A^1_x)^*\cong H^0(U, A^1)^* \qquad x \mapsto
[f \mapsto f(x)]
$$
which factors through the subspace
$$
T(x) := \suchthat{ F \in H^0(U), A^1)^*}{ F(c) = c\, \forall c\in \R}.
$$
The latter is clearly a torsor under the vector subspace
$$
H^0(U, \Lambda^1)^* \cong \suchthat{ F \in H^0(U, A^1)^*}{ F(c) = 0}
\subset H^0(U, A^1)^*.
$$

For example, if $\Delta$ is a 2-sphere and $A^1$ is $j_* \Aff$, where
$\Aff$ is the sheaf of affine functions on the complement of finitely
many points $P \subset \Delta$, then if $x \notin P$ then
$e: U \hookrightarrow e(U) \subset (x)$ is an open embedding, but if
the monodromy of the sheaf of affine functions on $\Delta \wo P$ is
non-trivial around $x \in P$ then $e: U \to e(U) \subset T(x)$ has
positive-dimensional fibres.

Now let $(\Delta, A^1)$ be an affine simplicial space.  For any face
$\sigma \subset \Delta$, there exists a neighborhood
$O \supset \sigma$ such that $e(O)$ is an open of a polyhedral complex
in $T(x)$, which is the same for any $x \in \sigma$, on which
superforms are defined as in Definition~\ref{sf-emb}.

If $\tau \supset \sigma$ are two faces of $\Delta$ then the inclusion
$\St(\tau) \subset \St(\sigma)$ induces a natural map $T(y) \to T(x)$
for any $y \in \sigma, x \in \tau$. Clearly, $e(\St(\tau))$ is mapped
to $e(\St(\sigma))$ under this map. We call the pullback of a
superform $\eta$ defined on an open set of $\St(\sigma)$ along this
map its \emph{restriction} to $U \cap \St(\tau)$.

\begin{defn}[Superforms on an affine simplicial space]
  \label{sf-deltax}
  Let $U \subset \Delta$ be an open subset of an affine simplicial
  space, and denote $\Sigma_U$ the collection of faces $\sigma$ such
  that $U \cap \mathring{\sigma} \neq \emptyset$.  A \emph{superform
    $\eta$} on an open subset $U \subset \Delta_X$ is a collection
  $(\eta_\sigma)_{\sigma \in \Sigma_U}$, where $\eta_\sigma$ is a germ
  of a superform on a neighbourhood of $U \cap \mathring{\sigma}$ in
  $\St(\sigma)$, such that whenever $\sigma \subset \tau$, the
  restriction of $\eta_\sigma$ to $U \cap \St(\tau)$ is $\eta_\tau$.
\end{defn}

The differentials $d'$ and $d''$, as well as maps $J, N, \bar N$ give
rise to the corresponding maps of sheaves of superforms on $\Delta$.

\begin{prop}
  \label{resolutions}
  \, \noindent
  \begin{enumerate}
  \item\label{a1-res}
    $\Ker \{ d' d'': \Aa^{0,0}_X \to \Aa^{0,1}_X \} \cong A^1 \otimes
    \R$;
  \item\label{ddc-symm} $N ( d'd'' f) = 0$ for all $f \in \Aa^{0,0}$;
  \item for all $p \geq 0$
    $\Ker \{ d'': \Aa^{p,0}_X \to \Aa^{p,1}_X \} \cong \Lambda^p_X
    \otimes \R$;
  \item\label{poinc}
    $\Im \{ d'': \Aa^{p,q} \to \Aa^{p,q+1} \} = \Ker \{ d'':
    \Aa^{p,q+1} \to \Aa^{p,q+2} \}$ for all $p \geq 0$;    
  \end{enumerate}
\end{prop}

\begin{proof}
  The first statement follows from the
  fact that $\Ker d'd''$ consists of linear combinations of affine
  functions. Indeed,
  $$
  d'd'' f = \sum_{i=1}^n \sum_{j=1}^n \dfrac{\partial^2 f}{\partial
    x_i \partial x_j} d'x_i \wedge d''x_j
  $$
  and the vanishing of all second derivatives of $f$ is equivalent to
  $f$ being affine.

  Further,
  $$
  Nd'd''\alpha = \sum_{i=1}^n \sum_{j=1}^n \left(\dfrac{\partial^2
      f}{\partial x_i \partial x_j} - \dfrac{\partial^2 f}{\partial x_j
      \partial x_i} \right) d''x_i \wedge d'' x_j = 0.
  $$

  The last two statements are proved in \cite[Seciton~5.2]{sus22}.
\end{proof}

In particular, $(\Aa^{p,q}_{\Delta}, d'')$ is a resolution of
$\Lambda^p$ and
$$
H^q(\Delta, \Lambda^p_\Delta) = \dfrac{ \suchthat{\alpha \in
    \Aa^{p,q}(\Delta)}{ d''\alpha = 0} }{ d'' \Aa^{p,q-1}(\Delta) }
$$
We will further denote $H^{p,q}(\Delta) = H^q(\Delta, \Lambda^p)$.

\begin{defn}
  If $(\Delta, A^1)$ and $(\tDelta, \tilde{A}^1)$ are two affine
  simplcial spaces then a continuous map $f: \Delta \to \tDelta$ is
  called a \emph{morphism of affine simplicial spaces} if $f^*
  \tilde{A}^1 \cong A^1$.
\end{defn}

\subsection{Positive supercurrents}

If $\dim V=n$ then an $(n,n)$-superform on $V$ is given by the
expression
$$
f\cdot d'x_1 \wedge d''x_1 \wedge \ldots \wedge d'x_n \wedge d''x_n,
$$
where $x_1, \ldots, x_n$ is a basis of $V$, ordered accordingly to the
chosen orientation. An $(n,n)$-superform is called positive if
$f \geq 0$.  Call a $(p,p)$-superform $\eta$ \emph{symmetric} if
$\eta = J\eta$.

\begin{defn}[Positive superforms on a real vector space]
  A symmetric $(p,p)$-superform $\eta$ is \emph{weakly positive} if
  $$
  \eta \wedge \alpha_1 \wedge J\alpha_1 \wedge \ldots \wedge
  \alpha_{n-p} \wedge J\alpha_{n-p}
  $$
  for any collection of forms 
  $\alpha_1, \ldots, \alpha_{n-p} \in \Aa^{1,0}(V)$.

  A symmetric $(p,p)$-superform $\eta$ is \emph{(strongly) positive}
  if
  $$
  \eta = \sum_j f_j \cdot \alpha_{1j} \wedge J\alpha_{1j} \wedge \ldots
  \wedge \alpha_{pj} \wedge J\alpha_{pj},
  $$
  where $f_j \geq 0$ and $\alpha_{ij} \in \Aa^{1,0}(V)$.
\end{defn}

These definitions naturally extend to superforms on affine simplicial
spaces.

\begin{defn}[Positive superforms on affine simplicial spaces]
  Let $\Delta$ be an affine simplicial space. A symmetric superform
  $\eta \in \Aa^{p,p}(\Delta)$ is \emph{weakly, resp. strongly
    positive}, if for any $x \in \Delta$, the germ $\eta_x$ is the
  pullback $e^* \tilde\eta_x$ of a germ of a weakly, resp. strongly
  positive form on a polyhedral complex in $T(x)$.
\end{defn}

\begin{defn}[Supercurrents]
  Let $U \subset \Delta$ be an open set. Let $\Dd^{p,q}(U)$ be the
  space of $(p,q)$-superforms on $U$ with compact support. We denote
  $\Dd'^{p,q}$ the space of continuous functions on $\Dd^{p,q}$ and
  call its elements \emph{currents} on $U$. The supercurrents form a sheaf
  $\Dd'^{p,q}$ on $\Delta$.  
\end{defn}

\begin{defn}[Positive supercurrents]
  Let $U \subset \Delta$ be an open set and let $T \in \Dd'^{p,q}$ be
  a supercurrent. We call $T$ \emph{weakly, resp. strongly, positive} if
  $T\eta \geq 0$ for any strongly, resp. weakly postive form
  $\eta \in \Dd^{p,q}(U)$.  
\end{defn}

\begin{fact}
  \label{reflexive}
  Let $U$ be an open subset of an affine simplicial space
  $\Delta$. Then the space $\Aa^{p,q}(U)$ is reflexive, that is, the
  topological dual of $\Dd'^{p,q}(U)$ is naturally isomorphic to
  $\Aa^{p,q}(U)$. 
\end{fact}

\begin{lemma}
  \label{pos-form}
  Let $(\Delta, A^1)$ be a compact affine simplicial space such that
  $A^1$ is constructible with respect to a simplicial complex
  strucutre that contains finitely many faces. Then there exists a
  positive symmetric $(1,1)$-superform on $\Delta$.
\end{lemma}

\begin{proof}
  Consider the finite covering of $\Delta$ by open stars of its
  vertices $\St(i)$ and select a positive form
  $\psi_i \in \Aa^{1,1}(\St(i))$. These forms then can be glued
  together to obtain a global $(1,1)$-form $\psi$ using a partition of
  unity subordinate to the cover $\set{\St(i)}$.
\end{proof}

\section{Affine Bott-Chern cohomology}
\label{bott-chern}

From now on we denote $\Delta$ a 2-sphere endowed with an affine
structure away from a finite set $P \subset \Delta$, and we assume
that the monodromy of the affine structure around any point $P$ is
unipotent. We denote $A^1 = j_* \Aff$, where $\Aff$ the sheaf of
affine functions on $\Delta \wo P$ and $j: \Delta \wo P \to \Delta$ is
the open embedding.

\subsection{Modification $\tDelta$ and Serre duality}

The sheaves of superforms that were just defined can be used to define
groups analogous to Dolbeault cohomology groups:
$$
H^{p,q}(\Delta) = \dfrac{\suchthat{ \eta \in \Aa^{p,q}(\Delta)
  }{d''\eta = 0}}{ d''\Aa^{p,q-1}(\Delta)}
$$
The proof of the main theorem that will be explained further in the
paper relies on certain properties of these cohomology groups. In
particular, one expects a ``Serre duality'':
$$
H^{p,q}(\Delta) \cong (H^{2-p,2-q}(\Delta))^*.
$$
However, we have $\dim H^{0,2} = 1$, but $\dim H^{0,2} = 0$: indeed,
$\dim T(p)=1$ for any affine structure singularity $p \in P$, and
therefore there are no non-zero germs of $(2,0)$-superforms at any
$p \in P$, and hence no global $d''$-closed $(2,0)$-forms.

In order to remedy this situation, consider the following modification
$\tDelta$ of $\Delta$ which is obtained by gluing simplices to
$\Delta$, a simplex for each $p \in P$. The idea of this construction
stems from \cite[Section~6.3]{sus22}: it is the modification to the
dual intersection complex of the central fibre $Y$ of a Kulikov
degeneration that corresponds to blowing up of a $(-1)$-curve
$C \subset Y_i$ that defines the simple affine structure singularity
at the vertex $i$. Here we effectively retell this construction
without a reference to a particular triangulation of $\Delta$ and
without the requirement that $\Delta$ arises as a dual intersection
complex.

We will explain the gluing near one point, the procedure is identical
near each point. Recall that by assumption the monodromy matrix of the
connection $\nabla$ associated to the affine structure on
$\Delta \wo P$ is unipotent, and therefore the affine structure in a
punctured neighbourhood of $i$ is obtained from a non-singular affine
structure on a neighbourhood of $i$ using a single ``shear'', as
described in, for example, \cite[Section~6.4]{ks06} or
\cite[Definition~8.3]{aet19}. Consider the following model of the germ
of such affine structure singularity in $\R^2$ that has monodromy
matrix
$$
\left( \begin{array}{cc}
         1 & 1 \\
         0 & 1 \\
       \end{array} \right)
$$
in the standard basis. Let
$$
\begin{array}{l}
W^+ = \suchthat{ (x,y) \in \R^2 }{ x > 0, y \geq 0}, \qquad  W^- =
  \suchthat{ (x,y) \in \R^2 }{ x > 0, y \leq 0},\\
  R = W^+ \cap W^-, \\
\end{array}
$$
and define the sheaf $A^1$ on $\R^2$ as follows: for any open $U$ that
does not intersect the ray $R$ let $A^1(U)$ be the affine functions on
$U$. For any open $U$ such that $U \cap R \neq \emptyset$ let the
sections of $A^1$ over $U$ be piece-wise affine functions $f$ that are
affine on $U \cap W^+$ and $U \cap W^-$ and such that the function
$$
f'(x,y) = \left\{\begin{array}{ll}
                   f(x+y, y) & (x,y) \in W^+\\
                   f(x, y) & (x,y) \in W^-\\
                 \end{array}\right.
$$               
is affine in the standard affine structure on $\R^2$.

Pick a point $i \in P$ and let $U$ be a neighbourhood of $i$, then
there exists a neighbourhood $U$ of $i$ and a continuous injective map
$\iota: U \to \R^2$ such that $A^1_\Delta|_U = \iota^* A^1_{\R^2}$. We
will slightly abuse notation denoting $x, y$ the pull-backs of the
standard coordinate functions to $U$.

Let $j \in \Delta$ be the point such that $x(j)=0, y(j) = 1$. we may
assume that there is no affine structure singularity at $j$,
multiplying $x$ and $y$ by a scalar if necessary. Glue a 2-simplex to
$\Delta$ so that one of its sides is $ij$ and the third vertex $e$ and
two other sides, $ie$ and $je$ are not glued. Let $\pi$ be a
linear projection from $\tDelta \to \Delta$ that $e$ maps to the
middle of the edge $ij$.

Let us now define the sheaf $A^1_{\tDelta}$. Firstly, we put it to be
isomorphic to $A^1_{\Delta}$ away from $ije$. For any neighbourhood
$U$ of $i$ or $j$ we put
$A^1_{\tDelta}(U) = \pi^*\iota^* A^1_{\Delta}(\iota(\pi(U)))$.

For any point $p$ in the interior of $ij$ and an open neighbourhood
$U$ of $p$ that doesn't meet $ie$ and $je$ we model $U$ on an
open neighbourhood of the origin of a fan in $\R^3$ that is generated
vectors
$$
e_1=(1,0,0),\ \ e_2=(0,1,0),\ \ e_3=(-1,-1,0),\ \ e_4=(0,0,1),\ \ e_5=(0,0,-1),
$$
with the 2-dimensional cones being spans of pairs of vectors
$$
\set{e_1, e_4}, \set{e_1, e_5}, 
\set{e_2, e_4}, \set{e_2, e_5}, 
\set{e_3, e_4}, \set{e_3, e_5} 
$$
(this fan is matroidal, see \cite[Definition~4.11]{jss19}).  For $p$
in the interior of $ie$, resp. $je$, and an open neighbourhood $U$ of
$p$ that is contained in the interior of $ije$ and $ie$, resp. $je$,
we let $A^1(U)$ to be the restrictions of affine functions $f$ on
$ije$ such that $f(i)=f(j)$. Finally, let $A^1$ have no sections in a
neighbourhood of $e$.

One can check that the projection $\pi: \tDelta \to \Delta$ is a
morphism of affine simplicial spaces, that is, $\pi^* A^1_{\Delta}
\subset A^1_{\tDelta}$.

\begin{prop}
  \label{serre-dual}
  The pairing
  $$
  H^q(\tDelta, \Lambda^p) \otimes H^{2-q}(\tDelta, \Lambda^{2-p}) \to
  \R, \qquad [\alpha] \otimes [\beta] \to \int \alpha \wedge \beta
  $$
  is non-degenerate.
\end{prop}

\begin{proof}
  The pairing is well-defined by Stokes theorem.
  
  The proof follows the strategy of \cite[Theorem~4.26]{jss19}: it
  suffices to show that the pairing is non-degenerate on local charts.
  The only non-trivial case to consider is the neighbourhoods of
  points in the interior of intervals $ij$ in the notation above, and
  the pairing is non-degenerate in this case by
  \cite[Proposition~4.19]{jss19}, since the sheaves of superforms are
  the same as ones defined on a matroidal fan.
\end{proof}

\subsection{Affine Harwey-Lawson criterion}
\label{harvey-lawson}

In complex geometry a K\"ahler current is a ``singular'' version of
the notion of a K\"ahler form, and existence of a K\"ahler current
characterizes manifolds that are bimeromorphic to K\"ahler manifoldds
(Fujiki class $\Cc$ manifolds).

When one considers a general affine simplicial space, it is not clear
what the correct notion of a ``top degree superform'' is, and
therefore it is not clear into which space of supercurrents the space
of superforms of given bidegree naturally embeds. In particular, it is
not clear what the dimension of a ``K\"ahler supercurrent'' should
be. Instead of attempting to build a general theory, we propose here
an \emph{ad hoc} definition of K\"ahler supercurrents on spaces
$\Delta$ and $\tDelta$, assuming simply that K\"ahler supercurrents
are continuous functionals on $(1,1)$-forms.

We fix a not necessarily closed symmetric positive $(1,1)$-superforms
$\psi_{\Delta}$ and $\psi_{\tDelta}$ on $\Delta$ and $\tDelta$, which
exist by Lemma~\ref{pos-form}.

\begin{defn}[K\"ahler supercurrents]
  Let $\psi$ be a positive $(1,1)$-superform on $\Delta$ or
  $\tDelta$. A \emph{K\"ahler supercurrent} is a symmetric positive
  $d''$-closed (and hence $d'$-closed) $(1,1)$-supercurrent $T$ such
  that $T - \psi$ is positive.
\end{defn}

\begin{prop} 
  \label{kahler-pushf}
  If $T \in \Dd'^{1,1}(\tDelta)$ is a K\"ahler supercurrent, then $\pi_*
  T$ is a K\"ahler supercurrent.
\end{prop}

\begin{proof}
  Straightforward.
\end{proof}

\begin{defn}[Weakly nef supercurrents]
  A suppercurrent $T \in \Dd'^{1,1}$ is called \emph{weakly nef} if
  it is a limit of positive $(n-1,n-1)$-forms.
\end{defn}

\begin{lemma}
  \label{make-kahler-curr}
  Let $\theta \in \Aa^{1,1}(\Delta)$. If $(\theta, \omega) \geq 0$ for
  every positive $(1,1)$-superform $\omega$ such that
  $d'd''\omega=0$. Then there exists a supercurrent
  $\chi \in \Dd'^{0,0}(\Delta_X)$ such that
  $\theta + d'd''\chi \geq 0$.
\end{lemma}

\begin{proof}
  Let
  $$
  \begin{array}{lll}
    P & = &  \suchthat{ \omega \in \Aa^{1,1}(\Delta_X)}{ \omega > 0}\\
    V & = & \suchthat{ \omega \in \Aa^{1,1}(\Delta_X)}{
            d'd''\omega = 0} \\
  \end{array}
  $$
  By assumption $\theta$ defines a functional on $V$ which is
  non-negative on the cone $P \cap V$, which is open in $V$.

  If $\theta$ reaches 0 on a non-zero vector in $P$ then since it is
  non-negative it must vanish on $V$. In this case $\theta$ defines a
  continuous functional on $d'd'' \Aa^{1,1} \subset \Aa^{2,2}$ which
  can be extended to a continuous functional $\chi$ on $\Aa^{2,2}$. In
  other words there exists a supercurrent
  $\chi \in \Dd'^{0,0}(\Delta_X)$ such that $\theta = -d'd''\chi$, so
  $\theta + d'd''\chi=0$.
 
  If $\theta$ is strictly positive on $P \cap V$, the hyperplane
  $\theta^\perp \subset V$ is disjoint from the convex
  cone $P \cap V$ and therefore by Hahn-Banach theorem the functional
  defined by $\theta$ on $\Aa^{1,1}$ can be extended to a
  functional $\tau$ that vanishes on $\theta^\perp \cap V$ and is
  strictly positive on $P$.

  Since $\theta^\perp$ is codimension 1 in $V$, $\theta$ and $\tau$
  must be proportional on $V$. Therefore for some positive constant
  $c$ the functional $\theta - c \tau$ vanishes on $V$, and by the
  same reasoning as above, $\theta - c \tau = -d'd'' \chi$ for some
  supercurrent $\chi \in \Dd'^{0,0}(\Delta_X)$. It follows that
  $\theta + d'd''\chi = c\tau \geq 0$.
  
\end{proof}

\begin{prop}
  \label{affine-hl}
  The following conditions are equivalent:
  \begin{enumerate}
  \item\label{cur} there exists a K\"ahler supercurrent on $\tDelta$;
  \item\label{bnd} every weakly nef $d''$-exact
    $(1,1)$-supercurrent is zero;
  \item\label{form} there exists a symmetric $(1,1)$-superform $\alpha$ such that
    $(\alpha, \omega) \geq (\psi, \omega)$ for any positive symmetric
    $(1,1)$-superform $\omega$ such that $d'd''\omega=0$.
  \end{enumerate}
\end{prop}

\begin{proof}
  (\ref{cur}) $\Rightarrow$ (\ref{bnd}) Let $\tau$ be the K\"ahler
  supercurrent, so $\tau = \theta + d'd''\chi \geq \psi$ for some $\theta$
  such that $N\theta=0, d''\theta=0$. Let
  $T =  d'' S$, $NS=0$, be a weak nef supercurrent such that $d'd'' T =0$, so
  $T = \lim \omega_n$ where $d'd'' \omega_n = 0, \omega \geq 0$. We
  have
  $$
  (\theta + d'd''\chi, \omega_n)  = (\theta, \omega_n) \geq (\psi, \omega_n)
  $$
  so $(\theta, T) \geq (\psi, T)$ in the limit. Furthermore,  
  $$
  (\theta, T) = (\theta, d'' S) = 0
  $$
  Therefore $(\psi, T) \leq 0$. Since $\psi > 0$, we have $T=0$.

  (\ref{bnd}) $\Rightarrow$ (\ref{form}) Define
  $$
  W = \suchthat{ \omega \in V \cap P }{ (\omega, \psi) = 1}
  $$
  and let $Y$ be the closure of $W$ in $\Dd'^{1,1}(\tDelta)$. By
  assumption, no element of $Y$ is $d''$-exact. Therefore, by
  Hahn-Banach theorem there exists $\eta \in \Aa^{1,1}(\Delta)$ such
  that $\eta|_{Y} > 0$, $d''\eta=0$. Since $Y$ is weakly compact,
  $\eta$ regarded as a functional achieves a minimum on $Y$:
  $(\eta, T) \geq C, \forall T \in Y$, and since $\eta|_Y > 0$,
  $C > 0$. In particular,
  $$
  \dfrac{(\eta,\omega)}{(\psi, \omega)} \geq C 
  $$
  since $\omega/(\omega, \psi) \in Y$. Taking $\alpha = \eta / C$, we
  obtain
  $$
  (\alpha, \omega) \geq (\psi, \omega).
  $$

  (\ref{form}) $\Rightarrow$ (\ref{cur}) If
  $\alpha \in \Aa^{1,1}(\Delta)$ is a $d''$-closed symmetric form
  such that $(\alpha, \psi) \geq (\omega, \psi)$ for any
  $d'd''$-closed superform $\omega$, then $\alpha - \psi$ satisfies the
  conditions of Lemma~\ref{make-kahler-curr}, and therefore exists a
  supercurrent $\chi$ such that $\alpha + d'd''\chi \geq \psi$, so
  $\alpha + d'd''\chi$ is a K\"ahler supercurrent. 
  
\end{proof}

\subsection{Monodromy morphism}

\begin{lemma}
  \label{proj-n}
  Let $V$ be a vector space, then for all $p,q \geq 0$
  $$
  \begin{array}{ll}
    (N \circ \bar N \circ N)\alpha = (q+1)N\alpha, & \textrm{ for all }
    \alpha \in \Aa^{1,q}(V), \\
    (\bar N \circ N \circ \bar N)\alpha = (p+1)
    \bar N \alpha, & \textrm{ for all } \alpha \in \Aa^{p,1}(V).\\
  \end{array}
  $$
  In particular, $\id - \frac{1}{q+1}\bar N \circ N$ is a projection
  on $\Ker N \subset \Aa^{1,q}$, and
  $\id - \frac{1}{p+1}N \circ \bar N$ is respectively a projection on
  $\Ker \bar N \subset \Aa^{p,1}$.
\end{lemma}

\begin{proof}
  Let $x_1, \ldots, x_n$ be coordinate. Let $\alpha$ be a
  $(p,q)$-superform given by the expression
  $$
  \alpha = \sum_{i,|J|=q} f_{i,J} d'x_I \wedge d''x_J.
  $$
  We have then
  \begin{align*}
    N \alpha & = \sum_{i=1}^n\sum_{|J|=q+1} f_{i,J}\ d''x_i \wedge
               d''x_J = \sum_{|J|=q+1} \left(\sum_{j \in J}\sgn(j,J)
               f_{j,J \wo \{j\}}\right) d''x_J,\\
    \bar N (N \alpha) & = \sum_{i=1}^n \sum_{|J|=q} \sgn(i, J \cup
                        \{i\}) \sum_{j \in J \cup \{i\}} \sgn(j,J \cup
                        \{i\}) f_{j,J \cup \{i\} \wo \{j\}}\ d'x_i
                        \wedge d''x_J, \\
  \end{align*}
  Applying the first formula to the second, we further get
  $$
    N (\bar N (N \alpha)) = \sum_{|J|=q+1} \sum_{k \in J} \sum_{j \in
      J} \sgn(k,J) \sgn(k,J) \sgn(j,J) f_{j,J \wo \{j\}}
    d''x_J,
  $$
  which proves the first claim. The statement for $\bar N$ is proved
  analogously.

  Now $\Im(\id - \bar N \circ N) = \Ker N$ since by the
  claim just proved
  $$
  N \circ (\id - \frac{1}{q+1}\bar N \circ N) = N - N = 0
  $$
  and since $\id - \frac{1}{q+1}\bar N \circ N = \id$ on $\Ker
  N$. Moreover,
  \begin{multline*}
    (\id - \frac{1}{q+1}\bar N \circ N) \circ (\id - \frac{1}{q+1}\bar
    N \circ N) = \id - \frac{1}{q+1}\bar N \circ N - \frac{1}{q+1}\bar
    N \circ N + \\
    + \frac{1}{(q+1)^2}\bar N \circ N \circ \bar N
    \circ N = \id - \frac{2}{q+1}\bar N \circ N \circ N +
    \frac{1}{q+1}\bar
    N \circ N = \id -  \frac{1}{q+1} \bar N \circ N.\\
  \end{multline*}
  The proof that the operator $\id - \frac{1}{p+1}N \circ \bar N$ is
  a projection on $\Ker \bar N \subset \Aa^{p,1}$ is analogous.
  
\end{proof}

\begin{lemma}
  \label{n-and-diff}  
  $$
   d'\circ N - N\circ d' =d'', \qquad d''\circ \bar  N - \bar N \circ d''=d'.
  $$

\end{lemma}

\begin{proof}
  Let $x_1, \ldots, x_n$ be coordinates let
  $$
  \alpha = \sum_{I,J} f_{IJ} d'x_I \wedge d''x_J
  $$
  be a $(p,q)$-superform. Then
  \begin{multline*}
    d' \alpha = \sum_{|I|=p+1, |J|=q} \sum_{i \in I} \sgn(i,I)
    \dfrac{\partial f_{I \wo
        \{i\}, J}}{\partial x_i} d'x_I \wedge d''x_J,\\
  \end{multline*}
  \begin{multline*}
    d'' \alpha = \sum_{|I|=p, |J|=q+1} (-1)^{p-1}\sum_{j \in J} \sgn(j,J)
    \dfrac{\partial f_{I, J \wo \{j\}}}{\partial x_i} d'x_I \wedge d''x_J,\\
  \end{multline*}
  \begin{multline*}
    N d' \alpha = \sum_{|I|=p, |J|=q+1} \sum_{j \in J}
    (-1)^{p}\sgn(j, I \cup \{j\}) \sgn(j,J) \cdot \\
    \cdot \sum_{i \in I \cup \{j\}} \sgn(i,I \cup \{j\})
    \dfrac{\partial f_{I \cup \{j\} \wo \{i\}, J \wo \{j\}}}{\partial
      x_i} d'x_I \wedge d''x_J,
  \end{multline*}
  \begin{multline*}
    d'N \alpha = \sum_{|I|=p, |J|=q+1} \sum_{i \in I} \sgn(i,I)
    \sum_{j \in J}
    (-1)^{p - 1}\sgn(j, I \cup  \{j\} \wo \{i\}) \sgn(j,J) \cdot \\
    \cdot \dfrac{\partial f_{I \cup \{j\} \wo \{i\}, J \wo \{j\}}}{\partial
      x_i} d'x_I \wedge d''x_J,
  \end{multline*}
  \begin{multline*}
     d'N\alpha - Nd' \alpha = \sum_{|I|=p, |J|=q+1} \sum_{j \in J}
    \sum_{i \in I} \biggl( (-1)^{p-1}\sgn(i,I)\sgn(j, I \cup \{j\} \wo \{i\}) \sgn(j, J) -  \\
    - (-1)^p\sgn(j, I \cup \{j\}) \sgn(i, I \cup
    \{j\}\sgn(j, J)  
    \biggr) \dfrac{\partial f_{I \cup \{j\}
        \wo \{I\} , J \wo \{j\}}}{\partial x_i} d'x_I \wedge d''x_J + \\
    - \sum_{j \in J}(-1)^p\sgn(j, I \cup \{j\})\sgn(j, I \cup
    \{j\})\sgn(j,J) \dfrac{\partial f_{I \cup \{j\} \wo \{j\} , J \wo
        \{j\}}}{\partial x_j}
    d'x_I \wedge d''x_J.\\
  \end{multline*}
  The first sum in the expression for $Nd'\alpha + d'N\alpha$ vanishes
  and the second sum coincides with the expression for $d''\alpha$,
  proving the first claim. Indeed, as one observes easily 
    $$
  \begin{array}{l}
    \sgn(i, I) \sgn(i,\{i,j\}) = \sgn(i, I \cup \{j\}), \\ 
    \sgn(j, I \cup \{j\} \wo \{i\}) \sgn(j,\{i,j\}) = \sgn(j, I \cup \{j\}), \\ 
  \end{array}
  $$
  and since $\sgn(i, \{i,j\}) = -\sgn(j,\{i,j\})$ by definition,
  \begin{multline*}
    (-1)^p\sgn(j, I \cup \{j\}) \sgn(i, I \cup \{j\}\sgn(j, J) -
    (-1)^{p-1}\sgn(i,I)\sgn(j, I \cup \{j\} \wo \{i\}) \sgn(j, J) = \\
    = (-1)^p\sgn(j, I \cup \{j\}) \sgn(i, I \cup \{j\}\sgn(j, J) + \\
    - (-1)^{p-1}(\sgn(i,\{i,j\})\sgn(i,I \cup
    \{j\}))(\sgn(j,\{j,i\})\sgn(j, I \cup \{j\})) \sgn(j, J) = 0 \\
  \end{multline*}

  Now, $\bar N = J \circ N \circ J$, so
  $$
  d'' \bar N \alpha - \bar N d'' \alpha  = (J \circ (d' \circ N - N \circ d') \circ J)(\alpha) = (J \circ d'' \circ J) \alpha
  = d' \alpha.
  $$

\end{proof}

\subsection{Resolution of $A^1$}

By Proposition~\ref{resolutions}(\ref{a1-res}) the sheaf $A^1$ is the
kernel of the morphism $d'd'': \Aa^{0,0} \to \Aa^{1,1}$. In fact, this
morphism can be extended to a resolution, which will be crucially used
in the arguments in the Section~\ref{main-sec}.

\begin{prop}
  \label{poincare-ddc}
  \,  \noindent

  \begin{enumerate}
  \item\label{poinc-n}
    $\Im \{ d'': \Aa^{p,q} \to \Aa^{p,q+1} \} = \Ker \{ d'':
    \Aa^{p,q+1} \to \Aa^{p,q+2} \}$ for all $p \geq 0$;
  \item\label{ddc-lemma}$d'd''(\Ker \{N: \Aa^{p,0} \to \Aa^{p-1,1})) = \Ker d'' \cap
    \Ker \{N: \Aa^{p,1}\to \Aa^{p-1,2}\}$.
  \item
    $d''(\Ker \{N: \Aa^{1,q} \to \Aa^{0,q+1})) = \Ker d'' \cap \Ker
    \{N: \Aa^{1,q+1}\to \Aa^{0,q+2}\}$ for all $q \geq 1$;
  \end{enumerate}
\end{prop}

\begin{proof}
  It follows from Definition~\ref{sf-deltax} that all four statements
  reduce to correspnoding statements about superforms on a vector
  space. In particular, the first statement reduces to
  \cite[Lemma~1.10]{lag12} or \cite[Theorem~2.16]{jss19}.

  For the second statement, the inclusion from left to right follows
  from Proposition~\ref{resolutions}(\ref{ddc-symm}) and from right to
  left from \cite[Lemma~1.13]{lag12}.

  In the last statement, the inclusion from left to right follows from
  Lemma~\ref{n-diff}. Let
  $\alpha \in \Aa^{1,q}, d''\alpha=N\alpha = 0$ for some $q \geq
  2$. Then by Proposition~\ref{resolutions}(\ref{poinc}) there exists
  $\beta \in \Aa^{1,q-1}$ such that $d''\beta = \alpha$. By
  Lemma~\ref{n-diff} again, $d''N\beta = 0$, so by
  Proposition~\ref{resolutions}(\ref{poinc}) there exists
  $\gamma \in \Aa^{0,q-1}$ such that $d''\gamma = \beta$. Consider
  $\beta + d''\bar N \gamma$: we have
  $$
  d''(\beta - \dfrac{1}{q-1} d'' \bar N \gamma) = d'' \beta = \alpha,
  $$
  so it is left to prove that this superform belongs to $\Ker
  N$. Indeed, applying Lemmas~\ref{proj-n} and \ref{n-and-diff} we get
  \begin{multline*}
    N(\beta - \dfrac{1}{q-1} d'' \bar N \gamma) = N\beta -
    \dfrac{1}{q-1}(Nd' \gamma + N \bar N d'' \gamma) = \\
    = N\beta - \dfrac{1}{q-1}( d'N \gamma - d''\gamma  + N \bar N N
    \beta) = N\beta - \dfrac{1}{q-1}(-N\beta + q \cdot N\beta) = N\beta
    - N\beta = 0.
  \end{multline*}
\end{proof}

It follows that the complex of sheaves 
\begin{equation}
  0 \to \Aa^{0,0} \xrightarrow{d'd''} \Aa^{1,1}_{Ker N}
  \xrightarrow{d''} \Aa^{1,2}_{\Ker N} \xrightarrow{d''} \Aa^{1,3}_{\Ker
    N}  \to \ldots
  \label{a1-reso}
\end{equation}
is the resolution of $A^1$ and one observes easily that its natural
inclusion into $\Cone(N)$
$$
0 \to \Aa^{0,0} \oplus \Aa^{0,1} \xrightarrow{ (d'', N + d'') } \Aa^{0,1}
\oplus \Aa^{1,1} \xrightarrow{ (d'',N + d'') } \Aa^{0,2} \oplus \Aa^{1,2}
\to \ldots
$$
is a quasi-isomorphism. In particular, we have a long exact sequence
of cohomology associated to the short exact sequence of sheaves
$0 \to \Lambda^0\cong \R \to A^1 \to \Lambda^1 \to 0$:
\begin{equation}
  \ldots \to H^{0,1}(\tDelta_X) \xrightarrow{N} H^{1,0}(\tDelta_X) \to
  H^{1,1}_{BC}(\tDelta) \xrightarrow{i} H^{1,1}(\tDelta_X) \to
  H^{0,2}(\tDelta_X) \to \ldots\label{eq:triangle}  
\end{equation}
where
$$
H^{1,1}_{BC}(\tDelta) = \dfrac{ \suchthat{\alpha \in
    (\Dd^{1,1}_{\Ker N})'(\tDelta)}{d''\alpha = 0}}{ d'd''
  \Dd'^{0,0}(\tDelta)},
$$
and the morphism $i$ on the level of cocycles is induced by the simple
inclusion of the space of $d''$-closed $(1,1)$-superforms in $\Ker N$
into the the space of $d''$-closed $(1,1)$-superforms.

By \cite[Lemma~1.10]{lag11} the resolution (\ref{eq:triangle}) can be
replaced with a similar one, where supercurrents are used:
\begin{equation}
  0 \to \Dd'^{0,0} \xrightarrow{d'd''} \Dd'^{1,1}_{\Ker N}
  \xrightarrow{d''} \Dd'^{1,2}_{\Ker N} \xrightarrow{d''} \Dd'^{1,3}_{\Ker N} \xrightarrow{d''} \ldots
  \label{a1-reso-curr}
\end{equation}
and so the group $H^{1,1}_{BC}(\tDelta)$ can be computed with supercurrents.

By Serre duality for superforms (Proposition~\ref{serre-dual}) the sequence
\begin{equation}
  \ldots \to (H^{0,2}(\tDelta))^* \xrightarrow{N} (H^{1,1}(\tDelta))^*
  \xrightarrow{p} (H^1_{BC}(\tDelta))^* \to (H^{1,2}(\tDelta))^* \to
  \ldots\label{eq:tri-dual}
\end{equation}
is exact, where
$$
(H^{1,1}_{BC}(\tDelta))^* = \dfrac{ \suchthat{\alpha \in
    \Aa^{1,1}_{\Ker N}(\tDelta) }{ d'd''\alpha = 0 }}{ \suchthat{\alpha
    \in \Aa^{1,1}_{\Ker N}(\tDelta)}{ \exists \beta \in \Aa^{1,0}_{\Ker
      N}(\tDelta), d''\beta=\alpha} }
$$

\section{Cohomology of sheaves $\Lambda^p$ on $\tilde{\Delta}$}
\label{main-sec}

\subsection{Map $Q$}

\begin{lemma}
  \label{conjugation}
  \begin{enumerate}
  \item The map $j_2: H^{2,0}(\tDelta) \to H^{0,2}(\tDelta)$ induced by $J$
    is an isomorphism;
  \item the map $N: H^{1,0}(\tDelta) \to H^{0,1}(\tDelta)$ is injective.
  \end{enumerate}

\end{lemma}

\begin{proof}
  Let $\omega \in \Aa^{2,0}(\tDelta_X), d''\omega=0$. Then
  $d'\omega=0$ for dimension reasons and $d''J\omega=0$. If there
  exists a superform $\alpha \in \Aa^{0,1}, d''\alpha = J\omega$ then
  $$
  d''(\omega \wedge \alpha) = d''\omega \wedge \alpha + \omega \wedge
  d''\alpha = \omega \wedge J\omega
  $$
  By Stokes theorem $\int_{\tDelta_X} d''(\omega \wedge \alpha) = 0$,
  but since $\omega \wedge J\omega$ is a positive superform,
  $\omega = 0$.  Therefore the map $j_2$ is injective.  Reasoning in a
  symmetric way, we get that the inverse of the map $j_2$ is injective
  too.
  
  Let $\alpha \in \Aa^{1,0}(\tDelta), d''\alpha=0$, so $\alpha$
  represents a class in $H^{1,0}(\tDelta)$. If $N\alpha=d''f$ then
  $\alpha = d'f$ and $d'd''f=0$. But since $\tDelta$ is compact, $f$
  mult be 0, so $\alpha=0$. 
\end{proof}

For the rest of this section we assume that the map
$N: H^0(\tDelta, \Lambda^1) \to H^1(\tDelta, \R)$ is an isomorphism then
by Poincar\'{e} duality the map
$$
H^{1,2}(\tDelta_X) \to H^{2,1}(\tDelta_X)
$$
induced by $J$ is an isomorphism, and the map
$$
H^1(\tDelta, A^1) \to H^1(\tDelta, \Lambda^1)
$$
is an injection. In view of Lemma~\ref{conjugation}
this is true in particular when $H^1(\tDelta, \R)=0$.

Denote
$$
\begin{array}{lll}
  Z^{1,1}_{d'd''} & = & \suchthat{ \omega \in \Aa^{1,1}_{\Ker
                        N}(\tDelta) }{ d'd''\omega=0}, \\
  Z^{1,1}_{d''} & = & \suchthat{ \omega \in \Aa^{1,1}(\tDelta) }{  d''\omega=0}.
\end{array}
$$

\begin{prop}
  \label{map-q}
  Assume that $N: H^{1,0}(\tDelta) \to H^{0,1}(\tDelta)$ is an
  isomorphism.  There exists a linear map
  $Q: Z^{1,1}_{d'd''} \to Z^{1,1}_{d''}$ such that
  \begin{enumerate}
  \item $Q(\omega) = \omega + Nd'\alpha$ for some
    $\alpha \in \Aa^{1,0}(\tDelta)$;
  \item $Q(d'' \gamma) = d''\gamma$ for all
    $\gamma \in \Aa^{1,0}(\tDelta)$.
  \end{enumerate}  
\end{prop}

\begin{proof}
  We will define $Q(\omega) = \omega + Nd'\alpha$ for the unique form
  $\alpha \in \Aa^{1,0}(\tDelta)$ such that $\omega + Nd'\alpha$ is
  $d''$-closed. This immediately implies that
  $Q(d''\gamma)=d''\gamma$. This map is also linear since if
  $Q(\omega_i)=\omega_i + Nd'\alpha_i, i=1,2$ then clearly
  $\omega_1+\omega_2 + Nd'(\alpha_1+\alpha_2)$ and
  $c\cdot \omega_1+ Nd'(c\cdot \alpha_1)$ are $d''$-closed for any
  constant $c \in \R$.

  First let us show that such form $\alpha$ exists for any
  $d'd''$-closed $\omega \in \Ker N$. Consider
  $d'\omega \in \Aa^{2,1}$, then since it is $d''$-closed by
  assumption, it defines a class in $H^{2,1}(\tDelta)$. The form
  $Jd'\omega$ is a $d''$-coboundary in $\Aa^{1,2}(\tDelta)$. Since by
  by assumption of the lemma $J$ induces an isomorphism
  $H^{1,0}(\tDelta) \to H^{0,1}(\tDelta)$, so by Fact~\ref{serre-dual},
  $J$ induces isomorphism $H^{2,1}(\tDelta) \to H^{1,2}(\tDelta)$,
  $d'\omega$ is a $d''$-coboundary in $\Aa^{2,1}(\tDelta)$ and
  therefore there exists $\beta \in \Aa^{2,0}$ such that
  $d''\beta = d'\omega$. Then by Lemma~\ref{n-and-diff}
  $$
  d''N\beta = Nd''\beta = Nd'\omega= d'N\omega-d''\omega=-d''\omega
  $$
  The form $J\beta$ is $d''$-closed for degree reasons and therefore
  defines a class in $H^{0,2}(\tDelta)$. The image of this class under
  $J$ is represented by a form
  $\tilde{\beta} \in \Aa^{2,0}, d''\tilde{\beta}=0$.  Since $J$
  induces an isomorphism $H^{0,2} \to H^{2,0}(\tDelta)$ by
  Lemma~\ref{conjugation}, $J\beta$ is cohomologous to
  $J\tilde{\beta}$, and so there exists
  $\alpha \in \Aa^{1,0}(\tDelta)$,
  $J \beta=J\tilde{\beta} + d''J\alpha$, so
  $\beta = \tilde{\beta} + d'\alpha$, and since
  $d''\beta=d''d'\alpha$,
  $$
  d''(\omega + Nd'\alpha) = d''\omega - d''\omega = 0.
  $$

  Let us prove that such form $\alpha \in \Aa^{1,0}(\tDelta)$ is
  unique. Suppose that there two forms $\alpha_1, \alpha_2$ such that
  both $\omega + Nd'\alpha_i$, $i=1,2$ are $d''$-closed. Then
  $$
  d''Nd'(\alpha_1 - \alpha_2) = Nd''d'(\alpha_1 - \alpha_2) = 0
  $$
  Since $N: \Aa^{2,1}(\tDelta) \to \Aa^{1,2}(\tDelta)$ is injective,
  $d'(\alpha_1 - \alpha_2)$ represents a class in $H^{2,0}$ which is
  mapped to a trivial class in $H^{0,2}$. By Lemma~\ref{conjugation},
  the class of $d'(\alpha_1-\alpha_2)$ is trivial in $H^{2,0}$, so is
  represented by 0, concluding the proof of uniqueness.
  
\end{proof}

\subsection{Proof of Theorem~A}

\begin{prop}
  \label{thm-a}
  There exists a K\"ahler supercurrent on $\Delta$.
\end{prop}

\begin{proof}
  By Propositions~\ref{affine-hl} and \ref{kahler-pushf} it suffices
  to show that every weakly nef $d''$-exact supercurrent
  $T \in \Ker N \subset \Dd'^{1,1}$ is zero.

  Clearly, such current $T$ is $d''$-closed and therefore represents a
  class in $H^{1,1}_{BC}(\tDelta)$. By
  Proposition~\ref{map-q}, the composition of natural morphisms
  $$
  H^{1,1}_{BC}(\tDelta) \to H^1(\tDelta, \Lambda^1_{\tDelta}) \to
  H^1(\tDelta, \Lambda^1_{\tDelta})^* \to H^{1,1}(\tDelta, \Lambda^1_{\tDelta})
  $$
  defines an isomorphism. The image of $[T]$ in
  $H^{1,1}_{BC}(\tDelta)$ under this morphism vanishes, since $T$ is
  $d'd''$-closed. Therefore $[T] = 0 \in H^{1,1}_{BC}(\tDelta)$. Let
  $\pi: \tDelta \to \Delta$ be the natural projection, then $\pi_* T$
  is still $d'd''$-exact and positive, in particular,
  $\pi_* T = d'd''f$ for some global convex function
  $f \in \Aa^{0,0}(\Delta)$. But the only convex functions on $\Delta$
  are the constant functions, therefore, $\pi_* T=0$ and $T=0$.
\end{proof}

\newpage
\footnotesize

\bibliography{kahler}

\vspace{2ex}

\noindent {\sc Department of Mathematics\\
  KU Leuven\\
  Celestijnenlaan 200B \\
  B-3001 Leuven (Heverlee)\\
  Belgium\\}
{\tt dsustretov.math@gmail.com\\}

\end{document}